\documentclass[12pt]{article}
\usepackage{amsmath,amssymb,amsthm,enumerate} 
\usepackage{tikz,pgfplots}

\usepackage{dsfont}
\usepackage{cite}

\def\nn{\nonumber}

\def\g{\gamma}

\def\vk{\varkappa}

\def\s{\sigma}
\def\la{\lambda}

\def\wt{\widetilde}
\def\ov{\overline}

\def\p{\partial}

\def\BC{{\mathbb C}}
\def\BR{{\mathbb R}}
\def\BN{{\mathbb N}}

\def\cld{{\mathcal D}}

\def\cld{{\mathcal D}}

\def\im{{\rm Im\ }}

\def\diag{\mathrm{diag}}
\def\sgn{{\rm sgn}}
\newcommand{\E}{\mathrm{e}}
\newcommand{\I}{\mathrm{i}}

\newtheorem{Pa}{Paper}[section]
\newtheorem{Tm}[Pa]{{\bf Theorem}}

\newtheorem{Cy}[Pa]{{\bf Corollary}}

\newtheorem{Ee}[Pa]{{\bf Example}}

\newtheorem{Pn}[Pa]{{\bf Proposition}}

\newenvironment{dedication}
        {\vspace{1ex}\begin{quotation}\begin{center}\begin{em}}    
        {\par\end{em}\end{center}\end{quotation}}

\title{GBDT and explicit solutions for  the \\  matrix coupled dispersionless
equations \\ (local and nonlocal cases)}

\author{Roman O. Popovych and Alexander Sakhnovich}

\date{}
\begin{document}
\maketitle

\begin{dedication} 
\end{dedication}

\begin{abstract}   We introduce matrix coupled (local and nonlocal) dispersionless
equations, construct wide classes of explicit multipole solutions, give explicit expressions 
for the corresponding Darboux and wave matrix valued functions  and consider their asymptotics in some
interesting cases. We consider the scalar cases of coupled, complex coupled and nonlocal
dispersionless equations as well.
\end{abstract}

{MSC(2010): 35B06, 37K40}

\vspace{0.2em}

{\bf Keywords:} matrix coupled dispersionless equation, matrix nonlocal dispersionless equation,
complex dispersionless equation, Darboux matrix,
transfer matrix function,
wave fuction, exlicit solution, asymptotics.  

\section{Introduction}\label{Intro}
\setcounter{equation}{0}
The { coupled dispersionless equations}  (real and complex) are integrable systems, which are actively studied since
the important works \cite{KaKo, KoO} (see, e.g., \cite{AbM3, Chen, FMO, Ha, KaKo2, Ku} and various references therein).
These equations are of independent interest and play also an essential role in the study of the short pulse equations
(see \cite{KaKo, FMO} and the references therein). We consider first the {\it matrix generalization of the coupled dispersionless equations}
(MCDE):
\begin{align}& \label{n1}
R_x=\frac{(-1)^p}{2}(VV_t+V_tV), \quad V_{tx}=\frac{1}{2}(VR+RV) \qquad \Big(R_x=\frac{\p}{\p x}R \Big),
\\ & \label{n2} R(x,t)=\diag\{\rho_1(x,t), \, \rho_2(x,t)\}, \quad V(x,t)= \begin{bmatrix}0 & v_1(x,t)\\ {v_2(x,t)} & 0 \end{bmatrix},
\end{align}
where diag stands for the block diagonal matrix, the blocks $\rho_k$ are $m_k\times m_k$ matrix functions $(k=1,2; \,\, m_k>0)$,
$v_1$ is an $m_1 \times m_2$ matrix function, $v_2$ is an $m_2 \times m_1$ matrix function, and $R$ and $V$ are
$m\times m$ matrix functions $(m:=m_1+m_2)$. Clearly, it suffices for $p$ in \eqref{n1} to take one of the two values $\{0,1\}$, that is, $p$ is either $0$ or $1$.

It is easy to see that system \eqref{n1} is equivalent to the compatibility condition 
\begin{align}& \label{n5}
 G_t-F_x+[G,F]=0, \qquad [G,F]:=GF-FG
\end{align}
of the following auxiliary linear systems:
\begin{align}& \label{n3}
w_x=G(x, t,\la)w;\quad G(x,t, \la):=\frac{\I}{4} \la j-\frac{\I}{2}j^{p+1}V(x,t), \quad j:=\begin{bmatrix} I_{m_1} &0 \\ 0 & -I_{m_2} \end{bmatrix}
\end{align}
(where $I_k$ is the $k \times k$ identity matrix), and
\begin{align}& \label{n4}
 w_t=F(x,t, \la)w, \quad F(x,t,\la):=\big(-\I  jR(x,t)+j^{p} V_t(x,t)\big)\big/ \la.
\end{align}
The complex  coupled dispersionless equations
\begin{align}& \label{1}
\rho_x+\frac{1}{2}\vk \big(|v|^2\big)_t=0, \quad v_{tx}=\rho v \quad (\rho=\ov {\rho}, \quad \vk = \pm 1)
\end{align}
appear, when we set in \eqref{n1}, \eqref{n3}, and \eqref{n4}: 
\begin{equation} \label{n6}
m_1=m_2=1, \quad \rho_1=\rho_2=\rho=\ov {\rho}, \quad v_1= v,  \quad v_2=\ov{v}, \quad p=(1+\vk)/2.
\end{equation}
Here $\ov {\rho(x,t)}$  is the complex conjugate of $\rho(x,t)$. We note that the generalized coupled dispersionless system in \cite{KaKo, KaKo2}  is more general than 
MCDE  \eqref{n1}. However, MCDE is more concrete and it is well known also (see, e.g., \cite{AbPrTr}) that the matrix and multicomponent
generalizations are of interest in applications.

For the case of MCDE, we introduce the GBDT-version of the B\"acklund--Darboux transformation and construct wide classes
of explicit solutions and corresponding explicit expressions for the Darboux  and wave matrix functions $w$. 
Various versions of B\"acklund--Darboux transformations and related
commutation methods are presented, for instance, in \cite{Ci, D, GeT, Gu, KoSaTe, Mar, MS, Schie} (see also the references therein).
For  generalized B\"acklund--Darboux (GBDT) approach see, for instance, \cite{ALS94, ALS-JMAA, SaSaR}.
The  Darboux matrices for the generalized coupled dispersionless systems
given in \cite{KaKo2} were constructed in \cite{Ha}   by an iterative procedure and for a special case of {\it diagonal generalized eigenvalues}.
GBDT allows to achieve an essential progress in this respect since neither the diagonal structure of the generalized eigenvalues nor
iterative procedure are required there.

Nonlocal nonlinear integrable equations 
have been actively  studied during the last years (see the important papers \cite{AbM3, Chen, GadA, GerIv, GuP} and 
numerous references therein), starting from the  article \cite{AbM0} on the nonlocal nonlinear Schr\"odinger equation.
The nonlocal (scalar) dispersionless equations were considered in \cite{AbM3, Chen}. Here we consider the
nonlocal case $R(-x)=-R(x)^*$, $V(-x)=V(x)^*$, that is,  
\begin{align}& \label{n7}
R(-x)=-R(x)^*, \quad v_1(x)=v(x), \quad v_2(x)=v(-x)^*.
\end{align}
We develop further the nonlocal results from \cite{MiSa}, introduce GBDT for the nonlocal equations \eqref{n1}, \eqref{n7}
and construct the corresponding explicit solutions and wave functions. The explicit construction of the wave functions is new
even for the local and nonlocal scalar dispersionless equations.

In Subsection \ref{SubAs} we consider also asymptotics of the Darboux matrix functions (Darboux matrices)
and, correspondingly, of the wave matrix functions (wave functions)
in the case of explicit solutions.

In the paper, $\BN$ denotes the set of natural numbers, $\BR$ denotes the real axis, $\BC$ stands for the complex plane, and
$\BC_+$ ($\BC_-$) stands for the open upper (lower) half-plane. The spectrum of a square matrix $A$ is denoted by $\s(A)$.
The notation diag means that the matrix is diagonal (or block diagonal).
\section{GBDT  for the matrix \\ coupled dispersionless  equations}\label{GBDT}
\setcounter{equation}{0}
\subsection{Preliminaries}
GBDT, which we consider here, is a particular case of the GBDT  introduced in \cite[Theorem 1.1]{ALS-JMAA}.
After fixing some $n \in \BN$, each GBDT for MCDE \eqref{n1} is determined by the initial system \eqref{n1} itself
and by {\it five parameter matrices  with complex-valued entries}:  three $n\times n$ invertible parameter matrices $A_1$, $A_2$  and $S(0,0)$
($\det A_i \not=0$, $\, i=1,2;$
$\det S(0,0)\not=0)$,
 and two
$n \times m$ parameter matrices $\Pi_1(0,0)$ and $\Pi_2(0,0)$ such that
 \begin{align}& \label{6}
A_1S(0,0)-S(0,0)A_2=\Pi_1(0,0)\Pi_2(0,0)^*.
\end{align}
Similar to \cite{ALS-JMAA}, we introduce coefficients $q_1(x,t)$, $q_0(x,t)$ and $Q_{-1}(x,t)$
via $G$ and $F$:
 \begin{align}& \label{7}
G=-\la q_1-q_0, \quad F=-\frac{1}{\la}Q_{-1}.
\end{align}
Hence, in view of \eqref{n3} and \eqref{n4} we have
 \begin{align} & \label{8}
 q_1=-\frac{\I}{4} j, \quad q_0=\frac{\I}{2}j^{p+1}V, \quad
Q_{-1}=\I  j R-j^p V_t.
\end{align}
If \eqref{n1} holds (and $V_t$ is continuous with respect to both variables combined), then the following linear differential systems are compatible and (jointly with the initial values
$S(0,0)$, $\Pi_1(0,0)$, and $\Pi_2(0,0)$) determine matrix functions $S(x,t)$, $\Pi_1(x,t)$, and $\Pi_2(x,t)$,
respectively:
\begin{align} \label{p1}& 
(\Pi_{1})_x= \sum _{i=0}^{1}A_{1}^{i} \Pi_{1}q_{i},
\quad (\Pi_{1})_t= A_{1}^{-1}
\Pi_{1}Q_{-1};
\\ & \label{p2}
(\Pi_{2})_x=- \sum _{i=0}^{1}(A_{2}^{*})^{i}
\Pi_{2}q_{i}^{*}, \quad (\Pi_{2})_t=-
(A_{2}^{*})^{-1} \Pi_{2}Q_{-1}^{*};
\\ & \label{s} 
S_{ x}= \Pi_{1}q_{1} \Pi _{2}^{*}, \quad   S_{ t}= - A_{1}^{-1}
\Pi_{1}Q_{-1} \Pi _{2}^{*} A_{2}^{-1}. \end{align}
Although the point $x=0$, $t=0$ is chosen above as the initial point, it is easy to see that
any other point may be chosen for this purpose as well.
Consider $S(x,t)$, $\Pi_1(x,t)$, and $\Pi_2(x,t)$ in some domain $D$, for instance, 
$$
D=\{(x,t):\, -\infty \leq a_1<x<a_2\leq \infty, \,\,  -\infty \leq b_1<t<b_2\leq \infty\},
$$
such that $R(x,t)$ and $V(x,t)$ of the form \eqref{n2}  are well defined in $D$ and satisfy \eqref{n1}, and such that $(0,0)\in D$.
Then $S(x,t)$, $\Pi_1(x,t)$, and $\Pi_2(x,t)$ are well defined and the identity
\begin{align}& \label{n8}
A_1S(x,t)-S(x,t)A_2=\Pi_1(x,t)\Pi_2(x,t)^*
\end{align}
follows from \eqref{6} and \eqref{8}--\eqref{s} \cite{ALS-JMAA}.

Introduce (in the points of invertibility of $S(x,t)$ in $D$) the matrix functions
\begin{equation} \label{n9}
w_{A}(x,t,\la)=I_{m}- \Pi_{2}(x,t)^{*}S(x,t)^{-1}(A_{1}- \la
I_{n})^{-1} \Pi_{1}(x,t).
\end{equation}
In view of \eqref{n8}, the matrix function $w_A(x,t, \la)$ is the so called transfer
matrix function in Lev Sakhnovich's form  (see \cite{SaSaR, SaL1, SaL3} and the references therein)
at each point $(x,t)$ of invertibility of $S(x,t)$.
Our next proposition refers to the particular case of \cite[(1.34)]{ALS-JMAA}.
\begin{Pn} \label{PnGF} Let $R$ and $V$ have the form \eqref{n2} and satisfy the MCDE  \eqref{n1}. 
Assume that $S(x,t)$, $\Pi_1(x,t)$ and $\Pi_2(x,t)$ satisfy \eqref{6} and \eqref{p1}--\eqref{s},
and that $w_A$ is given by \eqref{n9}.
Then, in the points of invertibility of $S(x,t)$ we have
\begin{align} & \label{n10}
\frac{\p}{\p x}w_{A}(x,t,z)=\wt G(x,t,z)w_{A}(x,t,z)-w_{A}(x,t,z)G(x,t,z), \\
 & \label{n11}
 \frac{\p}{\p t}w_{A}(x,t,z)=\wt F(x,t,z)w_{A}(x,t,z)-w_{A}(x,t,z)F(x,t,z),
\end{align}
where $G$ and $F$ are given by \eqref{7} and \eqref{8},
 \begin{align}& \label{7'}
\wt G=-\la \wt q_1-\wt q_0, \quad \wt F=-\frac{1}{\la}\wt Q_{-1}.
\\  & \label{8'}
 \wt q_1=q_1=-\frac{\I}{4} j, \quad \wt q_0=q_0-(q_1 X_0-X_0 q_1), 
 \\ \label{8!} &
\wt Q_{-1}=(I_m-X_{-1})Q_{-1}(I_m+Y_{-1}), 
\end{align}
$q_0$ and $Q_{-1}$ are defined by \eqref{8}, and $X_0$, $X_{-1}$ and $Y_{-1}$ are given by the formulas
\begin{equation}
 \label{8+}  X_0:=\Pi_2^*S^{-1}\Pi_1, \quad  X_{-1}:=\Pi_2^*S^{-1}A_1^{-1}\Pi_1, \quad Y_{-1}:=\Pi_2^*A_2^{-1}S^{-1}\Pi_1.
\end{equation}
\end{Pn}
According to \eqref{7} and \eqref{7'}, $G$ and $\wt G$ and $F$ and $\wt F$ depend on $\la$ in a similar way
(e.g., $G$ and $\wt G$ are polinomials of the first degree with respect to $\la$ and their leading coefficients coincide).

Compare \eqref{n9} and \eqref{8+} in order to see that $w_A(x,t,0)=I_m-X_{-1}(x,t)$. Hence, the equality \cite[(1.76)]{SaSaR} (at $z=\zeta=0$)
yields 
\begin{align}& \label{n12}
(I_m-X_{-1})(I_m+Y_{-1})=I_m.
\end{align}
Finally, according to \cite[(1.17)]{ALS-JMAA} the following useful relations are valid:
\begin{align}& \label{n13}
(\Pi_2^*S^{-1})_x=-\wt q_1 \Pi_2^*S^{-1}A_1-\wt q_0 \Pi_2^*S^{-1}, \quad (\Pi_2^*S^{-1})_t=-\wt Q_{-1} \Pi_2^*S^{-1}A_1^{-1}.
\end{align}
\subsection{Darboux matrix}
It follows from \eqref{7}, \eqref{8} and \eqref{7'}, \eqref{8'} that   $\wt q_0$ has the same form
as $q_0$, and so $\wt G$ has the same form as $G$. 
More precisely, in the expression for $\wt q_0$ we substitute only $\wt V$ instead of $V$, $\wt v_1$ instead of $v_1$ and $\wt v_2$ instead of $v_2$, 
that is
\begin{align}& \label{n14}
\wt q_0=\frac{\I}{2}j^{p+1}\wt V,
\end{align}
where
\begin{equation} \label{n15}
\wt{V}=\left[ \begin{array}{lr} 0 & \wt{v}_1
\\ \wt{v}_2 & 0 \end{array} \right]:=V+\frac{1}{2}j^p( X_0 -jX_0 j), \quad (X_0=\Pi_2^*S^{-1}\Pi_1).
\end{equation}
According to \eqref{7}, \eqref{8} and \eqref{7'}, the proposition below means that $\wt F$ has the same form as $F$ and $\wt Q_{-1}$ has the same form
as $Q_{-1}$.
\begin{Pn} \label{PnStrQ} Let the conditions of Proposition \ref{PnGF} hold.
Then, in the points of invertibility of $S$, we have
\begin{align}& \label{n16}
\wt Q_{-1}=\I  j \wt R-j^p \wt V_t.
\end{align}
where $\wt V$ is given by \eqref{n15} and
\begin{align}& \label{n17}
\wt R=\frac{1}{2 \I}( \wt Q_{-1}j+j \wt Q_{-1}), \quad \wt Q_{-1}=(I_m-X_{-1})Q_{-1}(I_m+Y_{-1}).
\end{align}
\end{Pn}
 \begin{proof} The second equality in \eqref{n17} coincides with \eqref{8!}. In view of the first equality in \eqref{n17} the block diagonal part of $\wt Q_{-1}$ equals
 $\I  j \wt R$ and in order to prove \eqref{n16} it remains to show that the block  antidiagonal part of $\wt Q_{-1}$
 equals $-j^p \wt V_t$, that is,
\begin{align}& \label{n18}
\frac{1}{2 }( \wt Q_{-1}-j \wt Q_{-1}j)=-j^p \wt V_t.
\end{align} 
First, let us find the derivative $(\Pi_2^*S^{-1}\Pi_1)_t$. Second equalities in \eqref{p1} and \eqref{n13} and the definition of $X_{-1}$ in \eqref{8+} yield
\begin{align} \nn
(\Pi_2^*S^{-1}\Pi_1)_t&=\Pi_2^*S^{-1} A_{1}^{-1}
\Pi_{1}Q_{-1}-\wt Q_{-1} \Pi_2^*S^{-1}A_1^{-1}\Pi_1
\\ & \label{n19}
=X_{-1}Q_{-1}-\wt Q_{-1}X_{-1}.
\end{align}
Using the second equality in \eqref{n17} and formula \eqref{n19} we derive
\begin{align} \label{n20}
(\Pi_2^*S^{-1}\Pi_1)_t
&=X_{-1}Q_{-1}-(I_m-X_{-1})Q_{-1}(I_m+Y_{-1})X_{-1}
\\ & \nn
=
X_{-1}Q_{-1}+(I_m-X_{-1})Q_{-1}(I_m+Y_{-1})(I_m-X_{-1})-\wt Q_{-1}.
\end{align}
By virtue of \eqref{n12}, we rewrite \eqref{n20} in the form
\begin{align} & \label{n21}
(\Pi_2^*S^{-1}\Pi_1)_t
=
X_{-1}Q_{-1}+(I_m-X_{-1})Q_{-1}-\wt Q_{-1}=Q_{-1}-\wt Q_{-1}.
\end{align}
Now, formulas \eqref{n15} and \eqref{n21} yield
\begin{equation} \label{n22} 
\wt{V}_t=V_t+\frac{1}{2}j^p( Q_{-1} -j Q_{-1} j)-\frac{1}{2}j^p(\wt Q_{-1} -j \wt Q_{-1} j).
\end{equation}
Hence, in view of the second equality in \eqref{8} we have
$$\wt{V}_t=-\frac{1}{2}j^p(\wt Q_{-1} -j \wt Q_{-1} j),$$ 
which implies \eqref{n18}. 
 \end{proof}

Propositions \ref{PnGF} and \ref{PnStrQ} lead us to the following theorem.
\begin{Tm}\label{Tm1} Let $R$ and $V$ have the form \eqref{n2}, let  $V_t(x,t)$ be a continuous function of $x$ and $t$ combined, and let $R$ and $V$
satisfy MCDE \eqref{n1} in $D$.  Assume that  three $n\times n$ parameter matrices $A_1$, $A_2$  and $S(0,0)$
$(\det A_i \not=0$, $\, i=1,2;$
$\det S(0,0)\not=0)$,
 and two
$n \times m$ parameter matrices $\Pi_1(0,0)$ and $\Pi_2(0,0)$ are given, and that the matrix identity \eqref{6} holds.

Then, $\Pi_1(x,t)$, $\Pi_2(x,t)$, $S(x,t)$ and $w(x,t,\la)$ $($where $w$ is the wave function, i.e., $w(x,t,\la)$  satisfies 
\eqref{n3}, \eqref{n4} and $\det w(0,0,\la)\not=0)$ are well defined in $D$. Moreover, in the points of invertibility of $S(x,t)$ in $D$, the matrix functions $\wt R$ and $\wt V$
given by \eqref{n15} and \eqref{n17}, respectively, have the form \eqref{n2} and satisfy \eqref{n1}, that is
 \begin{align}& \label{n1'}
\wt R_x=\frac{(-1)^p}{2}(\wt V \wt V_t+\wt V_t \wt V), \quad \wt V_{tx}=\frac{1}{2}(\wt V \wt R+\wt R \wt V),
\\ & \label{n2'} \wt R(x,t)=\diag\{\wt \rho_1(x,t), \, \wt \rho_2(x,t)\}, \quad \wt V(x,t)= \begin{bmatrix}0 & \wt v_1(x,t)\\ {\wt v_2(x,t)} & 0 \end{bmatrix}.
\end{align}
The wave function $\wt w$ $(\det \wt w(0,0,\la)\not=0)$, which corresponds to the transformed MCDE \eqref{n1'},  is given by the product $w_A\, w$:
\begin{align}& \label{n23}
\wt w(x,t,\la)=w_A(x,t,\la)w(x,t,\la); \quad   \wt w_x= \wt G \wt w, \quad   \wt w_t= \wt G \wt w;
\\ & \label{n24}
\wt G(x,t, \la)=\frac{\I}{4} \la j-\frac{\I}{2}j^{p+1}\wt V(x,t), 
\\ & \label{n25}
 \wt F(x,t,\la)=\big(-\I  j \wt R(x,t)+j^{p}\wt V_t(x,t)\big)\big/ \la .
\end{align}
 \end{Tm}
 \begin{proof}
 It follows from \cite{SaA-Comp} that $\Pi_1$, $\Pi_2$, $S$ and $w$ are well defined. Then, according to \eqref{n3}, \eqref{n4} and Proposition \ref{PnGF},
 $\wt w$ determined by the first equality in \eqref{n23} satisfies the second and third equalities in \eqref{n23}, where $\wt G$ and $\wt F$
 are given by \eqref{7'}. Moreover, the second and third equalities in
 \eqref{n23} imply that the compatibility condition 
 $$\wt G_t-\wt F_x+[\wt G, \wt F]=0$$ 
 holds.  Relations \eqref{7'}, \eqref{8'} and \eqref{n12} imply that \eqref{n24} holds. Relations \eqref{7'}  and \eqref{n16} yield \eqref{n25}.
 Moreover, according to the first equalities in \eqref{n15} and \eqref{n17}, the matrix functions $\wt R$ and $\wt V$ have the form \eqref{n2'}.
 
 Taking into account \eqref{n2'}, \eqref{n24} and \eqref{n25}, one can see that $\wt G$ and $\wt F$ have the same structure as $G$ and $F$, respectively.
 Thus, 
 the compatibility condition $\wt G_t-\wt F_x+[\wt G, \wt F]=0$ yields \eqref{n1'} in the same way as \eqref{n5} yields \eqref{n1}.
 \end{proof} 
 
 In particular, it is shown in Theorem \ref{Tm1} that the Darboux matrix, which transforms the wave function of the initial system into the
 wave function of the transformed system is given by the transfer function $w_A(x,t,\la)$.
 
 It is convenient to partition both $\Pi_1$ and $\Pi_2$ into $n\times m_1$ and $n\times m_2$ blocks:
  \begin{align}& \label{E0}
 \Pi_1=\begin{bmatrix}\Phi_1 & \Phi_2 \end{bmatrix}, \quad  \Pi_2=\begin{bmatrix}\Psi_1 & \Psi_2 \end{bmatrix}.
 \end{align} 
 The simplest cases where explicit solutions appear are the cases $V=0$ and $R=I_m$, $R=j$, $R=\I I_m$
 or $R=\I j$. For instance, when  $V=0$ and $R=I_m$ we  obtain $($in view of \eqref{8}--\eqref{p2}$)$ that
 \begin{align}& \label{E2}
\Phi_1(x,t)=\exp\big\{-\I\big((x/4)A_1-t A_1^{-1}\big)\big\}\Phi_1(0,0), 
\\ & \label{E3}
 \Phi_2(x,t)=\exp\big\{\I\big((x/4)A_1-t A_1^{-1}\big)\big\}\Phi_2(0,0),
 \\ & \label{E4}
 \Psi_1(x,t)=\exp\big\{-\I\big((x/4)A_2^*-t (A_2^*)^{-1}\big)\big\}\Psi_1(0,0),
  \\ & \label{E5}
 \Psi_2(x,t)=\exp\big\{\I\big((x/4)A_2^*-t (A_2^*)^{-1}\big)\big\}\Psi_2(0,0).
 \end{align}
 
 \begin{Ee}\label{Ee1} 
 Let us consider the case of trivial $V$ and constant diagonal matrix $R$ $(R=\cld)$
 with the entries $d_i$ $($or, written in the block form, blocks $\cld_1$ and $\cld_2)$ on the main diagonal$:$
 \begin{align}& \label{E1}
 V(x,t)\equiv 0, \quad R(x,t)\equiv \cld=\diag\{d_1,d_2, \ldots, d_m \}=\diag\{\cld_1, \cld_2\}, 
\\ & \label{E1'} 
  \cld_1=\diag\{d_1, \ldots, d_{m_1}\},
 \quad \cld_2=\diag\{d_{m_1+1}, \ldots, d_{m}\}.
 \end{align}
We set also 
 \begin{align}& \label{E6}
n=1, \quad A_1=a_1\in  \BC\setminus\{0\}, \quad A_2=a_2\in \BC\setminus\{0\} , \quad a_1\not=a_2.
 \end{align} 
 Then relations \eqref{8}--\eqref{p2}, \eqref{E1} and \eqref{E6} yield
 \begin{align}& \label{E2'}
\Phi_1(x,t)=\Phi_1(0,0)\exp\big\{-\I\big((x/4)a_1I_{m_1}-(t/a_1)\cld_1\big)\big\}, 
\\ & \label{E3'}
 \Phi_2(x,t)=\Phi_2(0,0)\exp\big\{\I\big((x/4)a_1I_{m_2}-(t/a_1)\cld_2\big)\big\},
 \\ & \label{E4'}
 \Psi_1(x,t)=\Psi_1(0,0)\exp\big\{-\I\big((x/4)\ov{a}_2I_{m_1}-(t/\ov{a}_2)\cld_1^*\big)\big\}, 
  \\ & \label{E5'}
 \Psi_2(x,t)=\Psi_2(0,0)\exp\big\{\I\big((x/4)\ov{a}_2I_{m_2}-(t/\ov{a}_2)\cld_2^*\big)\big\},
 \end{align} 
 where $\Phi_i$ and $\Psi_i$ are vector $($row$)$ functions.
  The function $S(x,t)$ may be recovered from \eqref{n8} and \eqref{E2'}--\eqref{E5'}$:$
 \begin{align}& \label{E7}
S(x,t)= (a_1-a_2)^{-1}(\Phi_1(x,t)\Psi_1(x,t)^*+\Phi_2(x,t)\Psi_2(x,t)^*),
\\ &  \label{E8}
\Phi_1(x,t)\Psi_1(x,t)^*
\\ \nn &
=\Phi_1(0,0)\exp\big\{-\I\big((x/4)(a_1-a_2)I_{m_1}-((t/a_1)-(t/a_2))\cld_1\big)\big\}\Psi_1(0,0)^*,
\\ &  \label{E9}
\Phi_2(x,t)\Psi_2(x,t)^*
\\ \nn &
=\Phi_2(0,0)\exp\big\{\I\big((x/4)(a_1-a_2)I_{m_2}-((t/a_1)-(t/a_2))\cld_2\big)\big\}\Psi_2(0,0)^*.
 \end{align} 
 Using \eqref{n15} and \eqref{n2'} we derive
 \begin{align}& \nn
\wt v_1(x,t)=\frac{1}{S(x,t)}\Psi_1(x,t)^* \Phi_2(x,t), \quad \wt v_2(x,t)=\frac{(-1)^p}{S(x,t)}\Psi_2(x,t)^* \Phi_1(x,t),
 \end{align} 
 where $\Phi_i$, $\Psi_i$ and $S$ are given explicitly in \eqref{E2'}--\eqref{E9}.
 Finally, from \eqref{8+}, \eqref{n17} and \eqref{n2'} we obtain
\begin{align} \nn
\wt \rho_1(x,t)=&\left(\begin{bmatrix}I_{m_1} & 0 \end{bmatrix}-\frac{1}{a_1S(x,t)}\Psi_1(x,t)^*\begin{bmatrix}\Phi_1(x,t) & \Phi_2(x,t) \end{bmatrix}\right)j\cld
\\ & \label{E10}\times
\left(\begin{bmatrix}I_{m_1} \\ 0 \end{bmatrix}+\frac{1}{a_2S(x,t)}\begin{bmatrix}\Psi_1(x,t)^* \\ \Psi_2(x,t)^* \end{bmatrix}\Phi_1(x,t)\right),
 \end{align} 
and we have a similar formula for $\wt \rho_2$ as well.
 Clearly, taking into account \eqref{n9}, \eqref{E0} and \eqref{E2'}--\eqref{E9} we have also an explicit formula for the Darboux matrix $w_A$.
\end{Ee} 
\subsection{Local matrix dispersionless equations  and \\ asymptotics of the Darboux matrix}\label{SubAs}
Let us set in \eqref{n2}
\begin{align}& \label{L1}
v_1(x,t)=v(x,t), \quad v_2(x,t)=v(x,t)^*; \quad R(x,t)=R(x,t)^*, \\
 \label{L1'} &
 {\mathrm{i.e.,}} \quad
 \rho_i(x,t)= \rho_i(x,t)^* \quad (i=1,2).
\end{align}
Then, MCDE \eqref{n1} takes the form of the  local matrix  dispersionless equation
\begin{align}& \label{L3}
v_{tx}=(\rho_1 v +v \rho_2)/2, \quad (\rho_1)_x=(-1)^p((vv^*)_t/2, \quad (\rho_2)_x=(-1)^p((v^*v)_t/2.
\end{align}
Put in \eqref{p1}--\eqref{n8},
\begin{align}& \label{L4}
A_1=A, \quad A_2=A^*, \quad \Pi_1(x,t)=\Pi(x,t), \quad \Pi_2(0,0)=-\I \Pi(0,0)j^{p+1}, \\
& \label{L5}
 S(0,0)=S(0,0)^*.
\end{align}
We will show that GBDT of the initial solutions of  \eqref{L3} into the transformed solutions of \eqref{L3} is determined 
by the triple of matrices  $\{A, S(0,0), \Pi(0,0)\}$, where $\det A\not=0$,
\begin{align}& \label{L7}
AS(0,0)-S(0,0)A^*=\I\Pi(0,0)j^{p+1}\Pi(0,0)^*,
\end{align}
and \eqref{L5} holds. 
According to \eqref{8} and \eqref{L1}, we have
\begin{align}& \label{L5'}
(Q_{-1}j^{p+1})^*=-Q_{-1}j^{p+1}.
\end{align}
It follows from \eqref{p1}, \eqref{p2} and from \eqref{L4}, \eqref{L5'} that
\begin{align}& \label{L6}
\Pi_2(x,t)=-\I \Pi(x,t)j^{p+1}.
\end{align}
Equations
\eqref{p1}--\eqref{n8} take the form
\begin{align}& \label{L8}
\Pi_x=A\Pi q_1+\Pi q_0, \quad \Pi_t=A^{-1}\Pi Q_{-1}, \\
& \label{L9}
 S_x=\Pi j^p \Pi^*/4, \quad
 S_t=A^{-1}\Pi \big(j^pR+(-1)^{p+1}\I jV_t\big)\Pi(A^*)^{-1},
 \\  & \label{L10}
 AS(x,t)-S(x,t)A^*=\I\Pi(x,t)j^{p+1}\Pi(x,t)^*.
\end{align}
Relations \eqref{L5} and \eqref{L9} yield
\begin{align}& \label{L6'}
S(x,t)=S(x,t)^*.
\end{align}
Next we show that 
\begin{align}& \label{L11}
\wt R(x,t)=\wt R(x,t)^*, \quad \wt V(x,t)=\wt V(x,t)^*,  
\end{align}
in the case considered in this subsection. In other words, if \eqref{L1} holds for the MCDE solutions,
then \eqref{L1} holds for the GBDT-transformed solutions as well. Indeed, relations \eqref{8!}, \eqref{8+},
\eqref{L4}, \eqref{L5'}, \eqref{L6}, and \eqref{L6'} imply that
\begin{align}& \label{L5+}
(\wt Q_{-1}j^{p+1})^*=-\wt Q_{-1}j^{p+1}.
\end{align}
From the first equalities in \eqref{n2} and \eqref{n17}, and from \eqref{L5+} we derive the first
equality in \eqref{L11}.  
Taking into account \eqref{L6}, we rewrite \eqref{n15} in the form
 \begin{align}& \label{L12}
\wt V=V+\frac{\I}{2}\left(j\Pi^*S^{-1}\Pi-\Pi^*S^{-1}\Pi j\right), \quad \left(V=\begin{bmatrix} 0 & v \\ v^* & 0 \end{bmatrix}\right),
\end{align}
and the second equality in \eqref{L11} follows.
Now, Theorem \ref{Tm1} yields the following corollary.
\begin{Cy}\label{LocMat} Let the $m_1 \times m_2$, $m_1\times m_1$ and $m_2\times m_2$ matrix functions $v(x,t)$, $\rho_1(x,t)$ and $\rho_2(x,t)$,
respectively, satisfy the local matrix dispersionless equation \eqref{L3} and  the equalities \eqref{L1'}, and let $v_t(x,t)$ 
be continuous in $D$.
Assume that the parameter matrices $A$ $(\det A\not=0)$, $S(0,0)=S(0,0)^*$ and $\Pi(0,0)$ satisfy \eqref{L7}.

Then, the matrix functions
$\wt v$, $\wt \rho_1$ and $\wt \rho_2$ given by \eqref{L12} and equalities
 \begin{align}& \label{L13}
\wt v=\begin{bmatrix} I_{m_1} & 0 \end{bmatrix}\wt V\begin{bmatrix}0 \\  I_{m_2}  \end{bmatrix}, \quad \begin{bmatrix} \wt \rho_1 & 0 \\
0 & \wt \rho_2
\end{bmatrix}=\frac{\I}{2\I}(\wt Q_{-1}j +j \wt Q_{-1}), \\
& \label{L14}
 \wt Q_{-1}=(I_m-\I j^{p+1} \Pi^* S^{-1}A^{-1}\Pi)(\I j R_t-j^p V_t)(I_m+\I j^{p+1}  \Pi^* (A^*)^{-1}S^{-1}\Pi),
\end{align}
where $S(x,t)$ and $\Pi(x,t)$ are determined by   \eqref{L8} and \eqref{L9}, satisfy the local matrix dispersionless equation and  the equalities
$\wt \rho_i(x,t)=\wt \rho_i(x,t)^* \quad (i=1,2)$.

The corresponding Darboux matrix $w_A$ takes the form
 \begin{align}& \label{L15}
w_A(x,t,\la)=I_m-\I j^{p+1}\Pi(x,t)^*S(x,t)^{-1}(A-\la I_n)^{-1}\Pi(x,t).
\end{align}
\end{Cy}

Further in the subsection, we consider the case
 \begin{align}& \label{L16}
R(x,t)\equiv R(t), \quad v(x,t)\equiv 0
\end{align}
and study the asymptotics of $v(x,t)$ and $w_A(x,t)$ when $x\to \infty$. {\it The asymptotics of $v(x,t)$ and $w_A(x,t)$ when $x\to -\infty$
can be studied in the same way.}

Formula \eqref{n23} for the fundamental solution $\wt w$ of the auxiliary systems (for the wave function) takes in this case
the form
\begin{align}& \label{L17}
\wt w(x,t,\la)=w_A(x,t,\la)\E^{\I \la x j/4}w(t,\la); \quad    w_t(t,\la)= \frac{1}{\I \la}jR(t)w(t,\la),
\end{align}
where $w_A$ is given by \eqref{L15}.  Hence, the asymptotics of the wave function with respect to $x$
is described by the asymptotics of the Darboux matrix $w_A$. Moreover, when we partition $\Pi(0,t)$ into
the $n\times m_1$ and $n \times m_2$ blocks $\Phi_1$ and $\Phi_2$, we have 
 \begin{align}& \label{L18}
\Pi(x,t)=\begin{bmatrix} \E^{xA/(4 \I)}\Phi_1(0,t) & \E^{-xA/(4 \I)}\Phi_2(0,t)\end{bmatrix}.
\end{align}
In view of \eqref{L9}, under the assumptions
 \begin{align}& \label{L19}
S(0,0)>0, \quad \sgn(t) j^p R(t) \geq 0
\end{align}
we have 
 \begin{align}& \label{L20}
S(0,t)>0; \quad S(x,t)>0 \,\, {\mathrm{for}} \,\, p=0, \,\, x\geq 0.
\end{align}
When $p=1$ and  \eqref{L19} holds, relations \eqref{L9} and \eqref{L10} yield
 \begin{align}& \label{L20'}
\big(\E^{\I  x A/4}S(x,t)\E^{-\I  x A^*/4}\big)'     \leq 0  , \quad \big(\E^{-\I  x A/4}S(x,t)\E^{\I  x A^*/4}\big)'\geq 0.
\end{align}
Since $S(0,t)>0$,  inequalities \eqref{L20'} imply that $\det S(x,t)\not=0$ and $\wt v(x,t)$ is well defined for
all $x\in \BR.$

We note that $w_A(4x,t,\la)\E^{\I \la x j}$ (for each fixed value $t$) is the fundamental solution
of the ``normalized" Dirac system
 \begin{align}& \label{L21}
y^{\prime}=\I(\la j +j^{p+1}\wt V_N(x))y, \quad \wt V_N= \begin{bmatrix} 0&-2\wt v(4x,t) \\ -2 \wt v(4x,t)^* & 0\end{bmatrix}.
\end{align}
Systems \eqref{L21} as the systems generated by the triples $\{A, S(0,t), \Pi(0,t)\}$ have been studied in a series
of papers (see \cite{FKRS, ALS-LAA16, SaSaR} and the references therein).
In particular, Weyl functions of the systems \eqref{L21} are rational, and  inverse problems to recover systems
from the rational Weyl functions have unique and explicit solutions.

Consider the case $p=0$. According to \cite[(3.13)]{ALS-WR} we have
\begin{align} \label{L22}&
\wt v(x,t)\in L^2_{m_1\times m_2}(\BR_+), \quad \wt v(x,t) \to 0 \quad {\mathrm{for}} \quad x \to \infty.
\end{align}
Without changing $\wt v$ and $w_A$ we may choose $n$, $A$, $S(0,t)>0$, and $\Pi(0,t)$ (see \cite{ALS-WR}) such that 
\begin{align}& \label{L23}
\s(A)\in (\BC_- \cup \BR), \quad {\mathrm{span}}\bigcup_{k=0}^{n-1}\im (A^k \Phi_1(0,t))=\BC^n,
\end{align}
where $\im$ stands for image and the second equality in \eqref{L22} means that the pair $\{A, \Phi_1(0,t)\}$ is controllable.
Then, we have the following asymptotic relation \cite[(3.28)]{ALS-WR}
\begin{align}\label{L24}&
w_{A}(x,t, \la )
=
\begin{bmatrix}
I_{m_1} &
0
\\
0 & \chi(t,\la)
\end{bmatrix} +o(1) \quad {\mathrm{for}} \quad x\to \infty, 
\\ & \label{L25} \chi(t,\la):=I_{m_2}+\I \Phi_2(0,t)^*\vk(t)(A-\la I_n )^{-1}\Phi_2(0,t),
\end{align}
where $\vk(t)=\lim_{x\to \infty}\Big(\E^{-\I x A } S(x,t)\E^{\I x A^{*} }\Big)^{-1}$, and this limit always exists.
\begin{Cy}
Let $p=0$ and assume that \eqref{L19} holds. Then $\wt v(x,t)$ does not have singularities when $x\geq 0$,
relations \eqref{L22} are valid and $($under a proper choice of  $n$, $A$, $S(0,t)>0$, and $\Pi(0,t))$ equality \eqref{L24} 
holds.
\end{Cy}
When \eqref{L19} holds (and so $S(0,t)>0$), one can combine \cite[Theorem 2.5]{ALS-LAA16} (see also the references
therein) and \cite[Theorem 3.7]{ALS-WR} in order to obtain the existence of the Jost solutions
\begin{align}& \label{L26}
\wt w_L'=(\I/4)\big(\la j-2 j \wt V(x,t)\big)\wt w_L; \quad \wt w_L(x,t,\la)= \E^{(\I /4)x\la j}(I_m+o(1))
\end{align}
for $\la \in \BR$ and $x \to \infty$. Moreover, from the above-mentioned theorems follows the expression for the
reflection coefficient $R_L$ of the form
\begin{align}& \nn
R_L(t,\la):=\begin{bmatrix}0 & I_{m_2}  \end{bmatrix} \wt w_L(0,t,\la)\begin{bmatrix}I_{m_1} \\ 0 \end{bmatrix}
\left(\begin{bmatrix}I_{m_1} & 0 \end{bmatrix} \wt w_L(0,t,\la)\begin{bmatrix}I_{m_1} \\ 0 \end{bmatrix}\right)^{-1}.
\end{align}
\begin{Cy} Let $p=0$ and assume that \eqref{L19} holds. Then the reflection coefficient $R_L$ of system
\eqref{L26} on the semi-axis $[0,\infty)$
is given by the formula
\begin{align}&\label{L27}
R_L(t,\la)=-\I \Phi_2(0,t)^*S(0,t)^{-1}(\la I_n-\theta)^{-1}\Phi_1(0,t), 
\end{align}
where $\theta=A-\I \Phi_1(0,t)\Phi_1(0,t)^*S(0,t)^{-1}$.
\end{Cy}

Consider the case $p=1$. We already showed that $\wt v(x,t)$ does not have singularities if \eqref{L19} holds.
Moreover, we have \cite[Corollary 3.6]{FKRS} 
$$\lim_{x\to \infty}\wt v(x,t)=0.$$ 
The asymptotics of $w_A$ may be
studied in a way similar to the case $p=0$ (see \cite{ALS-WR}) although the result is somewhat more complicated.

The  complex coupled dispersionless equation, which we consider in Section \ref{CCDE}, is a scalar subcase
of the local matrix dispersionless equation \eqref{L3}.
 \section{GBDT  for the nonlocal matrix  \\  dispersionless  equations}\label{NMDE}
\setcounter{equation}{0}
Recall that the nonlocal matrix dispersionless  equations are characterized  by the equalities \eqref{n7}. Equivalently,
the nonlocal  matrix dispersionless  equations (NMDE) are equations \eqref{n1}, where $R$ and $V$ have the form
\begin{align}& \label{n26}
 R(x,t)=\diag\{\rho_1(x,t), \, \rho_2(x,t)\}=-R(-x,t)^*, \\
& \label{n27} 
  V(x,t)= \begin{bmatrix}0 & v(x,t)\\ {v(-x,t)}^* & 0 \end{bmatrix} \quad (\mathrm{i.e.,} \quad V(x,t)=V(-x,t)^*).
\end{align}
In the nonlocal case we assume \eqref{n26} and \eqref{n27} instead of \eqref{n2} and (similar to the
subsection \ref{SubAs}) determine GBDT by 3 parameter matrices. However, these matrices
satisfy somewhat  different relations. Namely, we set
\begin{align}& \label{n28}
A_1=A \quad(\det A\not=0), \quad A_2=-A^*, \quad \Pi_1(x,t)=\Pi(x,t), \\
& \label{n28'}
 \Pi_2(0,0)=-\I \Pi (0,0) j^p, \quad S(0,0)=-S(0,0)^* \quad (\det S(0,0)\not=0).
\end{align}
 so that the identity \eqref{6} takes the form
\begin{align}& \label{n29}
AS(0,0)+S(0,0)A^*=\I\Pi(0,0)j^p \Pi(0,0)^* .
\end{align} 
It easily follows from \eqref{p1}--\eqref{s} that \eqref{n28} and \eqref{n28'} yield
\begin{align}& \label{n31}
\Pi_2(x,t)\equiv -\I \Pi(-x,t)j^p, \quad S(x,t)\equiv -S(-x,t)^*.
\end{align}  
Thus, the identity \eqref{n8} takes the form 
\begin{align}& \label{n32}
AS(x,t)+S(x,t)A^*=\I\Pi(x,t)j^p \Pi(-x,t)^* .
\end{align}  
In view of  \eqref{n31}, we have $X_0=\I j^p\Pi(-x,t)^*S(x,t)^{-1}\Pi(x,t)$ and formula \eqref{n15} takes the form
\begin{align} & \label{n33}
\wt{V}(x,t)=V+\frac{\I}{2}\Big( \Pi(-x,t)^*S(x,t)^{-1}\Pi(x,t) 
-j\Pi(-x,t)^*S(x,t)^{-1}\Pi(x,t) j\Big).
\end{align} 
Let us again partition $\Pi$ into two blocks: $\Pi=\begin{bmatrix}\Phi_1 & \Phi_2 \end{bmatrix}$, 
where $\Phi_1$ is an $m\times m_1$ matrix function. Now, \eqref{n27} and \eqref{n33} imply that
\begin{align} & \label{n34}
\wt{V}(x,t)=\left[ \begin{array}{lr} 0 & \wt{v}(x,t)
\\ \wt{v}(-x,t)^* & 0 \end{array} \right]=\wt V(-x,t)^*, \\
& \label{n35}
 \wt{v}(x,t)=v(x,t)+\I \Phi_1(-x,t)^*S(x,t)^{-1}\Phi_2(x,t).
\end{align} 

Next, we show that  $\wt R$ given by \eqref{n17} satisfies (under the assumptions of this section) the nonlocal requirement
\begin{align}& \label{n26'}
\wt R(x,t)=\diag\{\wt \rho_1(x,t), \, \wt \rho_2(x,t)\}=-\wt R(-x,t)^*.
\end{align}
Indeed, in view of  the relations \eqref{8+}, \eqref{n28},  and \eqref{n31}, we have
\begin{align}& \label{n36}
X_{-1}(x,t)=\I j^p\Pi(-x,t)^*S(x,t)^{-1}A^{-1}\Pi(x,t), \\
& \label{n36'}
 Y_{-1}(x,t)=-\I j^p\Pi(-x,t)^*(A^{-1})^*S(x,t)^{-1}\Pi(x,t),
\\ & \label{n37}
X_{-1}(-x,t)^*=-j^pY_{-1}(x,t)j^p, \quad Y_{-1}(-x,t)^*=-j^pX_{-1}(x,t)j^p
\end{align}
From the last equality in \eqref{8} and the formulas \eqref{n26} and \eqref{n27} we derive
\begin{align}& \label{n38}
Q_{-1}(-x,t)^*=j^pQ_{-1}(x,t)j^p.
\end{align}
Formulas \eqref{8!}, \eqref{n37} and \eqref{n38} imply that
\begin{align}& \label{n39}
\wt Q_{-1}(-x,t)^*=j^p\wt Q_{-1}(x,t)j^p.
\end{align}
Finally, the first equality in \eqref{n17} and formula \eqref{n39} yield \eqref{n26'}.

Rewriting \eqref{p1}, \eqref{s}, and \eqref{n9} under assumptions of this section, we obtain
\begin{align} \label{p1+}& 
\Pi_x= \sum _{i=0}^{1}A^{i} \Pi q_{i},
\quad \Pi_t= A^{-1}
\Pi Q_{-1}; \quad S_{ x}(x,t)= \I \Pi(x,t) q_{1} j^p \Pi(-x,t)^{*},
\\ & \label{s+} 
 S_{ t}= \I A^{-1}
\Pi(x,t)Q_{-1}(x,t)j^p \Pi(-x,t)^{*} (A^{*})^{-1},
\\ & \label{w}
w_{A}(x,t,\la)=I_{m}- \I j^p\Pi(-x,t)^{*}S(x,t)^{-1}(A- \la
I_{n})^{-1} \Pi(x,t).
\end{align}

Recall that in this section we assume that the relations \eqref{n28} and \eqref{n28'} hold (in particular,
GBDT is determined by the triple $\{A, S(0,0), \Pi(0,0)\}$). Now, we can rewrite Theorem \ref{Tm1} for the NMDE case.
\begin{Tm}\label{Tm2} Let $R$ and $V$ have the forms \eqref{n26} and \eqref{n27}, respectively, let $V_t(x,t)$ be continuous in $D$, and let $R$ and $V$
satisfy \eqref{n1} in $D$.  Assume that  two $n\times n$ parameter matrices $A$  and $S(0,0)$
 and one
$n \times m$ parameter matrix $\Pi(0,0)$  are given, and that  \eqref{n29} holds.

Then, $\Pi(x,t)$,  $S(x,t)$ and $w(x,t,\la)$ $($where $w$ is the wave function, i.e., $w(x,t,\la)$  satisfies 
\eqref{n3}, \eqref{n4} and $\det w(0,0,\la)\not=0)$ are well defined in $D$. 

Moreover, in the points of invertibility of $S(x,t)$ in $D$, the matrix function $\wt R$ 
given by \eqref{n17} $($where $X_{-1}$ and $Y_{-1}$ have the forms \eqref{n36} and \eqref{n36'}$)$ and the matrix function $\wt V$ given by \eqref{n33} satisfy NCDE,
that is, the relations \eqref{n26'} and \eqref{n34} are valid and the equations
 \begin{align}& \label{n1+}
\wt R_x=\frac{(-1)^p}{2}(\wt V \wt V_t+\wt V_t \wt V), \quad \wt V_{tx}=\frac{1}{2}(\wt V \wt R+\wt R \wt V)
\end{align}
are satisfied.

The wave function $\wt w$ $(\det \wt w(0,0,\la)\not=0)$, which corresponds to the transformed NCDE \eqref{n1+},  is given by the product $w_A\, w$,
where the Darboux matrix $w_A$ has the form \eqref{w}. In other words, the relations \eqref{n23}--\eqref{n25} are valid.
 \end{Tm}
 \section{GBDT  for the complex \\ coupled dispersionless  equations}\label{CCDE}
\setcounter{equation}{0}
Recall that in order to obtain the complex  coupled dispersionless  equations (CCDE)
\begin{align}& \label{1'}
\rho_x+\frac{1}{2}\vk \big(|v|^2\big)_t=0, \quad v_{tx}=\rho v \quad (\rho=\ov {\rho}, \quad \vk = \pm 1),
\end{align}
we set in the MCDE (see \eqref{n1} and \eqref{n2}) the equalities \eqref{n6}. In particular, since $m_1=m_2=1$, the functions $\rho$ and $v$ are scalar functions
and we rewrite \eqref{8} and \eqref{n2} in the form
\begin{align}
& \label{8!!}
 q_1=-\frac{\I}{4} j, \quad q_0=\frac{\I}{2}j^{p+1}V, \quad
Q_{-1}=\I  \rho j -j^p V_t;
\\ \label{n40} &
j=\begin{bmatrix}1 &0 \\ 0 & -1 \end{bmatrix}, \quad R(x,t)=\rho(x,t)I_2, \quad  V(x,t)=\begin{bmatrix}0 & v(x,t)\\ \ov{v(x,t)} & 0 \end{bmatrix}.
\end{align}
In view of \eqref{n40}, the coefficients $q_1,\, q_0, \, Q_{1}$ given by \eqref{8!!} have the property
 \begin{align}& \label{10}
q_k^*= -  j^{p+1}q_kj^{p+1}            \quad (k=1,0), \quad Q_{-1}^*=-j^{p+1}Q_{-1}j^{p+1}.
\end{align}
In order to construct GBDT for the CCDE equations \eqref{1'}, we set in the GBDT for MCDE in Section \ref{GBDT} the equalities
 \begin{align}& \label{n41}
A_1=A_2^*=A, \quad \Pi_1(x,t)=\Pi(x,t), \quad \Pi_2(0,0)=-\I \Pi(0,0)j^{p+1}, 
\end{align}
and $S(0,0)=S(0,0)^*$.  It means that GBDT is determined by 3 parameter matrices:
$$A \quad (\det A \not=0), \quad S(0,0)=S(0,0)^* \quad (\det S(0,0) \not=0), \quad
{\mathrm{and}} \quad \Pi(0,0).$$
 In view of \eqref{n41}, we rewrite \eqref{p1} in the form
\begin{align}& \label{9}
\Pi_x=A\Pi q_1+\Pi q_0, \quad \Pi_t=A^{-1}\Pi Q_{-1}.
\end{align}
Taking into account \eqref{p2} and \eqref{10}--\eqref{9}, we see that 
$$ \Pi_2(x,t)=-\I \Pi(x,t)j^{p+1}.$$
Hence, the identity \eqref{n8} takes the form
 \begin{align}& \label{12}
AS(x,t)-S(x,t)A^*=\I\Pi(x,t)j^{p+1}\Pi(x,t)^*.
\end{align}
The matrix function $S(x,t)=S(x,t)^*$ is determined now  by $S(0,0)$ and
the equations 
 \begin{align}& \label{11}
S_x=\I \Pi q_1 j^{p+1}\Pi^*, \quad S_t=-\I A^{-1}\Pi Q_{-1}j^{p+1}\Pi^*(A^*)^{-1},
\end{align}
which follow from \eqref{s}. 

The GBDT-transformed solution $\{\wt R, \, \wt V\} $, Darboux matrix  $w_A$ and wave function $\wt w$
are expressed via $A$, $\Pi(x,t)$ and $S(x,t)$.  Let us show that $\wt R$ and $\wt V$ have
the form \eqref{n40}:
\begin{align}
 \label{n42} &
\wt R(x,t)=\wt \rho(x,t)I_2  \quad (\wt \rho=\ov {\wt \rho}), \quad \wt V(x,t)=\begin{bmatrix}0 &\wt v(x,t)\\ \ov{\wt v(x,t)} & 0 \end{bmatrix},
\end{align}
and so $\wt \rho$ and $\wt v$ satisfy CCDE. Indeed, $X_0$, $X_{-1}$ and $Y_{-1}$ (given by \eqref{8+}) take now the form
\begin{align}&
 \label{8!!+}  X_0=\I j^{p+1}\Pi^*S^{-1}\Pi,  \\
 \label{8!!++}  &
 X_{-1}=\I j^{p+1}\Pi^*S^{-1}A^{-1}\Pi, \quad
 Y_{-1}=\I j^{p+1}\Pi^*(A^*)^{-1}S^{-1}\Pi.
\end{align}
Thus, we rewrite \eqref{n15} as
\begin{equation} \label{n43}
\wt{V}=V+\frac{\I}{2}( j \Pi^*S^{-1}\Pi -\Pi^*S^{-1}\Pi j).
\end{equation}
According to \eqref{n43}, $\wt V$ has the form \eqref{n42}, where
\begin{equation} \label{n44}
\wt{v}=v+\I \Phi_1^*S^{-1}\Phi_2 \qquad \big(\begin{bmatrix} \Phi_1 & \Phi_2 \end{bmatrix}:=\Pi\big).
\end{equation}
According to \eqref{8!}, \eqref{n12} and \eqref{8!!} we have
\begin{align}& \label{27}
{\mathrm{tr}}(\wt Q_{-1})={\mathrm{tr}}(Q_{-1})=0,
\end{align}
where ${\mathrm{tr}}$ stands for trace.
In view of \eqref{n17}  and \eqref{27}, the first equality in \eqref{n42} holds. 

It remains to prove that $\wt \rho=\ov {\wt \rho}$. From \eqref{8!!++} we see that 
$$X_{-1}^*=-j^{p+1}Y_{-1}j^{p+1}.$$
 Hence, the last equality in \eqref{10} and the equalities in
\eqref{n17} imply that $\wt Q_{-1}^*=-j^{p+1}\wt Q_{-1}j^{p+1}$ and so $\wt R=\wt R^*$.
That is, we have
\begin{align}& \label{29}
\wt \rho =\I\begin{bmatrix} 0 & 1 \end{bmatrix}\wt Q_{-1}\begin{bmatrix} 0 \\ 1 \end{bmatrix}=\ov{\wt \rho} \qquad \big(\wt Q_{-1}=(I_m-X_{-1})Q_{-1}(I_m+Y_{-1})\big),
\end{align}
which finishes the proof of \eqref{n42}. We obtained the following corollary of Theorem \ref{Tm1}.
\begin{Cy} Let $v_t(x,t)$ be continuous in $D$, and
let the functions $\rho$ and $v$ satisfy CCDE \eqref{1'} in $D$.  Assume that  two $n\times n$ parameter matrices $A$  and $S(0,0)=S(0,0)^*$
 and one $n \times 2$ parameter matrix $\Pi(0)$  are given, and that  the  relations
  \begin{equation} \label{n45}
AS(0,0)-S(0,0)A^*=\I\Pi(0,0)j^{p+1}\Pi(0,0)^* \quad (\det A\not=0, \quad \det S(0,0)\not=0)
\end{equation}
hold. Introduce $\Pi(x,t)$ and $S(x,t)$ using \eqref{9}, \eqref{11} $\big($and \eqref{8!!}, \eqref{n40}, where $p=(1+\vk)/2\big)$.

Then, in the points of invertibility of $S(x,t)$ in $D$, the functions $\wt \rho$ 
$($given by \eqref{29} and \eqref{8!!++}$)$ and $\wt v$ $($given by \eqref{n44}$)$  satisfy CCDE $:$
\begin{align}& \label{1+}
\wt \rho_x+\frac{1}{2}\vk \big(|\wt v|^2\big)_t=0, \quad \wt v_{tx}= \wt \rho \, \wt v \quad (\wt \rho=\ov {\wt \rho}).
\end{align}
Moreover, a wave function $w(x,t,\la)$  $($where $\det w(0,0,\la)\not=0)$ is well defined in $D$ via \eqref{n40} and auxiliary systems
\begin{align}& \label{2}
w_x=G(x, t,\la)w,\quad G(x,t, \la):=\frac{\I}{4} \la j-\frac{\I}{2}j^{p+1}V(x,t);
\\ & \label{3}
 w_t=F(x,t, \la)w, \quad F(x,t,\la):=\big(-\I \rho j+j^{p} V_t(x,t)\big)\big/ \la .
\end{align}
The wave function $\wt w$ $(\det \, \wt w(0,0,\la)\not=0)$, which corresponds to the transformed CCDE \eqref{1+},  is given by the product $\wt w=w_A\, w$,
where the Darboux matrix $w_A$ has the form 
 \begin{align}& \label{13}
w_A(x,t, \la)=I_2-\I j^{p+1}\Pi(x,t)^*S(x,t)^{-1}(A-\la I_n)^{-1}\Pi(x,t).
\end{align}
\end{Cy}
\begin{Ee}\label{Ee2} In order to  present an example of the solution of  CCDE \eqref{1'},  we set in Example \ref{Ee1}
$($in accordance with 
\eqref{n6} and \eqref{n41}$)$
$$m_1=m_2=1, \quad a_1=a, \quad a_2=\ov{a}, \quad \cld_1=\cld_2=d=\ov{d}, $$ 
$\Psi_1(0,0)=-\I \Phi_1(0,0),$ and 
$\Psi_2(0,0)=(-1)^p\I \Phi_1(0,0)$. For simplicity of notations, we put $\Phi_i(0,0)=c_i$. In view of \eqref{E7}--\eqref{E9}
we have
 \begin{align}\nn
S(x,t)=&\frac{\I}{a-\ov{a}}\Big(|c_1|^2\exp\big\{-\I\big((a-\ov{a})(x/4)-d((t/a)-(t/\ov{a}))\big)\big\}
\\ & \label{E11}
+(-1)^p|c_2|^2\exp\big\{\I\big((a-\ov{a})(x/4)-d((t/a)-(t/\ov{a}))\big)\big\},
 \end{align} 
and $S(x,t)=\ov{S(x,t)}$.  The formula for $\wt v_1$ in Example \ref{Ee1} takes the form
  \begin{align} & \label{E12}
\wt v(x,t)=\frac{\I\ov{c}_1 c_2}{S(x,t)}\exp\big\{\I\big((a+\ov{a})(x/4)-d((t/a)+(t/\ov{a}))\big)\big\}.
 \end{align} 
Finally, formula \eqref{E10}  takes the form
 \begin{align}\nn  
\wt \rho(x,t)=&d\Big(1-\I |c_1|^2\exp\big\{-\I\big((a-\ov{a})(x/4)-d((t/a)-(t/\ov{a}))\big)\big\}\big/(aS(x,t))\Big)
\\ & \nn \times
\Big(1+\I |c_1|^2\exp\big\{-\I\big((a-\ov{a})(x/4)-d((t/a)-(t/\ov{a}))\big)\big\}\big/(\ov{a}S(x,t))\Big)
\\ & \label{E13}
+(-1)^p d |\wt v(x,t)^2|\big/ |a|^2.
\end{align} 
\end{Ee}

\section{Coupled dispersionless  equations}\label{Li}
\setcounter{equation}{0}
Similarly to Section \ref{CCDE}, we consider here the case of scalar function $\rho$
(and scalar $v_1$ and $v_2$). Setting in \eqref{n1}
 \begin{align}& \label{n46}
m_1=m_2=1, \quad p=1, \quad \rho_1=\rho_2=\rho, \quad v_1=r, \quad v_2=s,
\end{align}
we rewrite \eqref{n1} in the form
 \begin{align}& \label{n47}
\rho_x+\frac{1}{2}(rs_t+r_ts)=0, \quad r_{tx}=\rho r, \quad s_{tx}=\rho s,
\end{align}
which is equivalent, for instance, to \cite[(1.2)]{Li} (see also the references therein). The following corollary of Theorem \ref{Tm1} is valid.
\begin{Cy}\label{CyLi} Let the conditions of Theorem \ref{Tm1} hold and assume additionally that
 \begin{align}& \label{n48}
m_1=m_2=1,  \quad \rho_1=\rho_2=\rho.
\end{align}
Then,
$\wt v_1$ and $\wt v_2$ given by  \eqref{n15}, and $\wt \rho$ given by the formula
\begin{align}& \label{n50}
\wt \rho =\I\begin{bmatrix} 0 & 1 \end{bmatrix}\wt Q_{-1}\begin{bmatrix} 0 \\ 1 \end{bmatrix} \qquad \big(\wt Q_{-1}=(I_2-X_{-1})Q_{-1}(I_2+Y_{-1})\big)
\end{align}
satisfy equations \eqref{n47}, that is,
 \begin{align}& \label{n49}
\wt \rho_x=\frac{(-1)^p}{2}\big(\wt v_1(\wt v_2)_t+(\wt v_1)_t\wt v_2\big), \quad (\wt v_1)_{tx}=\rho \wt v_1, \quad (\wt v_2)_{tx}=\rho \wt v_2.
\end{align}
\end{Cy}
\begin{proof} Taking into account \eqref{n1}, \eqref{n2} and \eqref{n17}, the only fact which we need to prove
is that $\wt R$ has the form $\wt \rho I_2$, that is, that  $\wt \rho_1=\wt \rho_2$. Similar to the calculation
in the previous section, relations \eqref{8} and \eqref{n12} yield
\begin{align}& \label{n51}
{\mathrm{tr}}(\wt Q_{-1})={\mathrm{tr}}(Q_{-1})=0,
\end{align}
and the equality $\wt \rho_1=\wt \rho_2$ follows from \eqref{n17} and \eqref{n51}.
\end{proof}

For the nonlocal situation
\begin{align}& \label{n52}
\ov{\rho(x,t)}=-\rho(-x,t), \quad  v(x,t):=v_1(x,t)=\ov{v_2(-x,t)},
\end{align}
equations \eqref{n49} have  the form
\begin{align}& \label{n53}
 \rho_x(x,t)=\frac{(-1)^p}{2}\big( v(x,t)\ov{v_t(-x,t)}+ v_t(x,t)\ov{v(-x,t)}\big), \\
 & \label{n54}
  v_{tx}(x,t)=\rho(x,t)  v(x,t).
\end{align}
In other words, under conditions \eqref{n48} and \eqref{n52} system \eqref{n1} is equivalent to the system \eqref{n53}, \eqref{n54}.
(Note that $\rho$ and $v$ are scalar functions.) 

Assume further that the relations \eqref{n28}--\eqref{n29} hold (in particular,
GBDT is determined by the triple $\{A, S(0,0), \Pi(0,0)\}$). Below, we formulate a corollary of Theorem \ref{Tm2}.
\begin{Cy}\label{CyLi2} Let  $v_t$ be continuous in $D$, assume that $\ov{\rho(x,t)}=-\rho(-x,t)$,
and let $v$ and $\rho$ satisfy \eqref{n53}, \eqref{n54} in $D$.

Then, $\wt v$ given by \eqref{n35} and $\wt \rho$ given by \eqref{n50}
satisfy \eqref{n53}, \eqref{n54}  $($as well as  $v$ and $\rho)$, that is, the following equalities hold$:$
\begin{align}& \label{n53'}
\wt \rho_x(x,t)=\frac{(-1)^p}{2}\big( \wt v(x,t)\ov{\wt v_t(-x,t)}+ \wt v_t(x,t)\ov{\wt v(-x,t)}\big), \\
 & \label{n54'}
 \wt v_{tx}(x,t)=\wt \rho(x,t) \wt v(x,t).
\end{align}
We also have $\ov{\wt \rho(x,t)}=-\wt \rho(-x,t)$.
\end{Cy}
\begin{proof} Similar to the Corollary \ref{CyLi} we need only to prove that $\wt \rho_1=\wt \rho_2$
(after which we may use Theorem \ref{Tm2}).  The equality $\wt \rho_1=\wt \rho_2$ follows from \eqref{n17} and \eqref{n51}.
\end{proof}
\section{Examples and figures}\label{ExFig}
In these examples we construct explicit solutions of the nonlocal equations \eqref{n53'}, \eqref{n54'}. 
We set
\begin{align} & \label{E14}
v(x,t)\equiv 0, \quad R(x,t) \equiv \I I_2, \,\, {\mathrm{i.e.}} \,\, \rho(x,t) \equiv \I.
\end{align} 
We put
\begin{align} & \label{E15}
n=2, \quad A=\begin{bmatrix} a & 1\\ 0 & a \end{bmatrix} \quad (a+\ov{a}\not=0),
\end{align} 
which corresponds  \cite{ALS-WR} to the simplest case of the Weyl function $($reflection coefficient$)$ with a pole of the order
more than one $($so called multipole case$)$. For the literature on the multipole cases see, for instance, \cite{Ol, Schie} and the references
therein. 

Using notations $\Phi_i(0,0)=C_i=\begin{bmatrix}c_{i1} \\ c_{i2}\end{bmatrix}$, 
we  $($similar to the deduction of \eqref{E2}$)$  obtain
 \begin{align}\nn
\Phi_1(x,t)&=\exp\big\{-\big(\I (x/4)A+t A^{-1}\big)\big\}C_1
\\ &  \label{E16}
=\exp\big\{-\big((\I/4) a x+(t/a)\big)\big\}\big(I_2-(\I/4)xA_0+(t/a^2)A_0\big)C_1,
 \end{align}
where $A_0=\begin{bmatrix}0 & 1 \\ 0 & 0 \end{bmatrix}$.
In the same way as \eqref{E16}, we derive
\begin{align}&  \label{E18}
\Phi_2(x,t)
=\exp\big\{(\I/4) a x+(t/a)\big\}\big(I_2+(\I/4)xA_0-(t/a^2)A_0\big)C_2.
 \end{align} 
Relations \eqref{E16}--\eqref{E18} provide an explicit expression for 
\begin{align} & \label{E19}
K(x,t)=\{K_{i\ell}(x,t)\}_{i,\ell=1}^2=\I\big(\Phi_1(x,t)\Phi_1(-x,t)^*+(-1)^p\Phi_2(x,t)\Phi_2(-x,t)^*\big).
\end{align}
Next, using \eqref{n32} we easily express $S(x,t)=\{S_{i\ell}(x,t)\}_{i,\ell=1}^2$ in terms of  $K(x,t):$
\begin{align} & \label{E20}
S_{22}(x,t)=K_{22}(x,t)/(a+\ov{a}), 
\\ & \label{E20'} S_{12}(x,t)=\big(K_{12}(x,t)-S_{22}(x,t)\big)/(a+\ov{a}), 
\\ & \label{E21}
S_{21}(x,t)=\big(K_{21}(x,t)-S_{22}(x,t)\big)/(a+\ov{a}), 
\\ & \label{E22}
S_{11}(x,t)=\big(K_{11}(x,t)-S_{12}(x,t)-S_{21}(x,t)\big)/(a+\ov{a}).
\end{align}
Finally, from \eqref{n35}, \eqref{n50} and \eqref{E14} it follows that
\begin{align}  \label{n35'}
 \wt{v}(x,t)=&\I \Phi_1(-x,t)^*S(x,t)^{-1}\Phi_2(x,t),
\\ \nn
\wt \rho(x,t)=&\I\big(1-\I(-1)^p\Phi_2(-x,t)^*S(x,t)^{-1}A^{-1}\Phi_2(x,t)\big)
\\ \nn & \times
\big(1-\I(-1)^p\Phi_2(-x,t)^*(A^*)^{-1}S(x,t)^{-1}\Phi_2(x,t)\big)
\\ \nn & 
+\I (-1)^p\Phi_2(-x,t)^*S(x,t)^{-1}A^{-1}\Phi_1(x,t)
\\ \label{E23} & \times
 \Phi_1(-x,t)^*(A^*)^{-1}S(x,t)^{-1}\Phi_2(x,t).
\end{align} 
Recall that $($according to Corollary \ref{CyLi2}$)$  $\wt v$ and $\wt \rho $ satisfy \eqref{n53'}, \eqref{n54'}.

The   fundamental solution $($wave function$)$ $w(x,t,\la)$ of the initial systems \eqref{n3} and \eqref{n4},
where $V=0$ and $R=\I I_2$ is given by the
formula 
$$w(x,t,\la)=\exp\big\{\big((\I/4)\la x+(t/\la)\big)j\big\}.$$
In view of \eqref{E15}--\eqref{E22} we have explicit formulas
for the Darboux matrix $w_{A}(x,t,\la)=I_{2}- \I j^p\Pi(-x,t)^{*}S(x,t)^{-1}(A- \la
I_{2})^{-1} \Pi(x,t)$. Thus the wave function $w_A(x,t,\la)w(x,t,\la)$ of the transformed
system with $\wt v $ and $\wt \rho$ given by \eqref{n35'} and \eqref{E23}, respectively, is also expressed
explicitly.

Let us consider several explicit formulas in greater detail. In the following we set $a_1=\Re(a)$ and $a_2=\Im(a)$.

\begin{figure}
\centering
\includegraphics[scale=1.0]{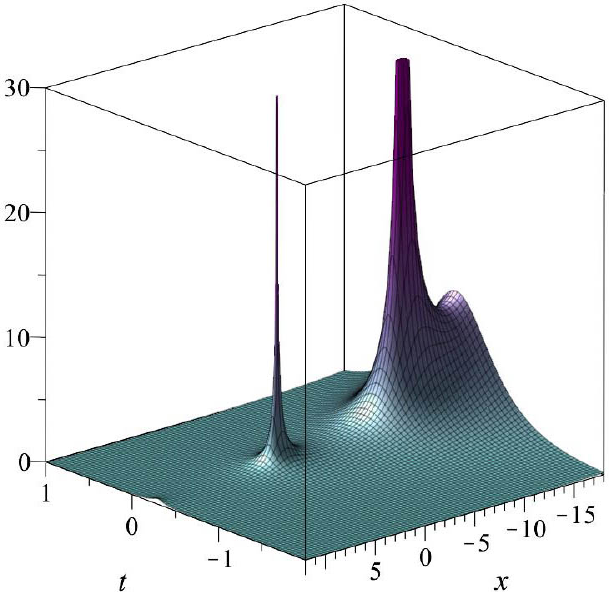}\qquad 
\includegraphics[scale=1.0]{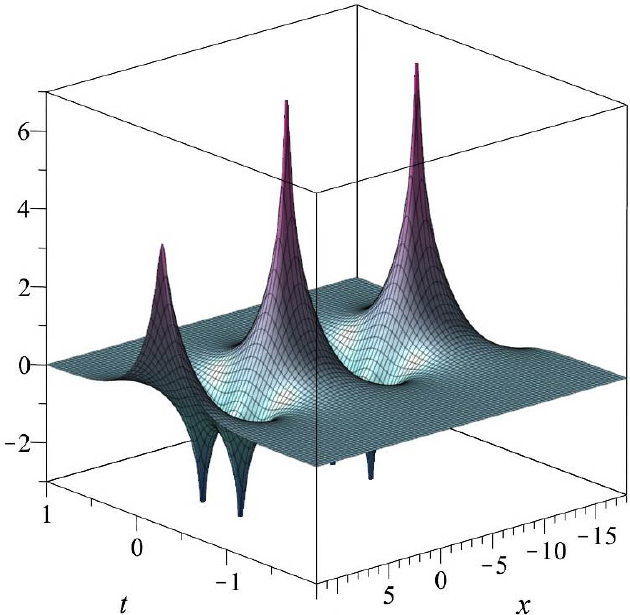}    
\caption{$|\wt v|$ (left) and $\ln|\wt\rho|$ (right). }\label{fig1}
\end{figure}

{\it Case 1.} The simplest case is the case where $c_{12}=c_{21}=0$, that is, $C_1=\begin{bmatrix} c_{11}\\ 0 \end{bmatrix}$ and 
$C_2=\begin{bmatrix} 0 \\ c_{22} \end{bmatrix}$. Here, relations \eqref{E16}--\eqref{E23} after some calculations
yield:
\begin{gather*}
\wt v=\frac {4a_1^2\bar c_{11}c_{22}\E^{-a_2(x+4\I t/|a|^2)/2}} {4a_1^2|c_{11}|^2\E^{-\I a_1(x-4\I t/|a|^2)/2}+(-1)^p|c_{22}|^2\E^{\I a_1(x-4\I t/|a|^2)/2}},
\\[1ex]
\wt\rho=\I -
\frac
{32\I (-1)^pa_1^4|c_{11}|^2|c_{22}|^2/|a|^2}
{\big( 4a_1^2|c_{11}|^2\E^{-\I a_1(x-4\I t/|a|^2)/2} + (-1)^p|c_{22}|^2\E^{\I a_1(x-4\I t/|a|^2)/2} \big)^2}.
\end{gather*}
In particular, for 
$$p=1, \quad a=\frac{1}{2}+\frac{1}{3}\I, \quad c_{11}=1+2\I,  \quad c_{22}=4+3\I,$$
the behaviour of $|\wt v|$ and $\ln |\wt \rho|$ is shown on Figure \ref{fig1}.

\newpage

\begin{figure}[h!]
\centering
\includegraphics[scale=1.0]{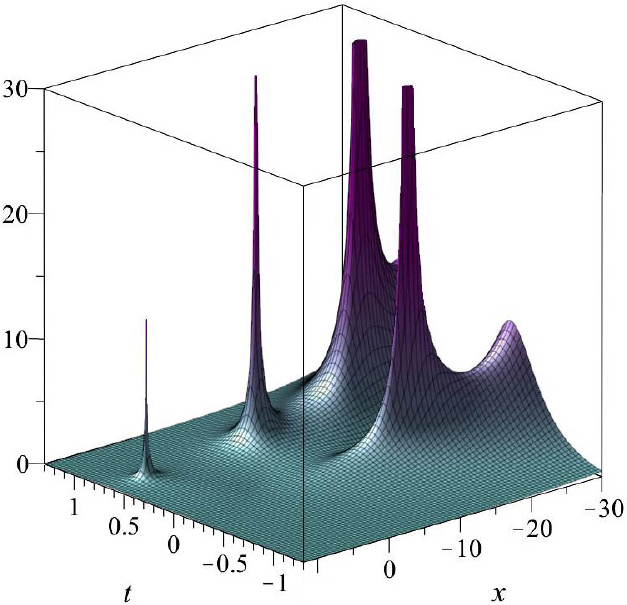}\quad\ 
\includegraphics[scale=1.0]{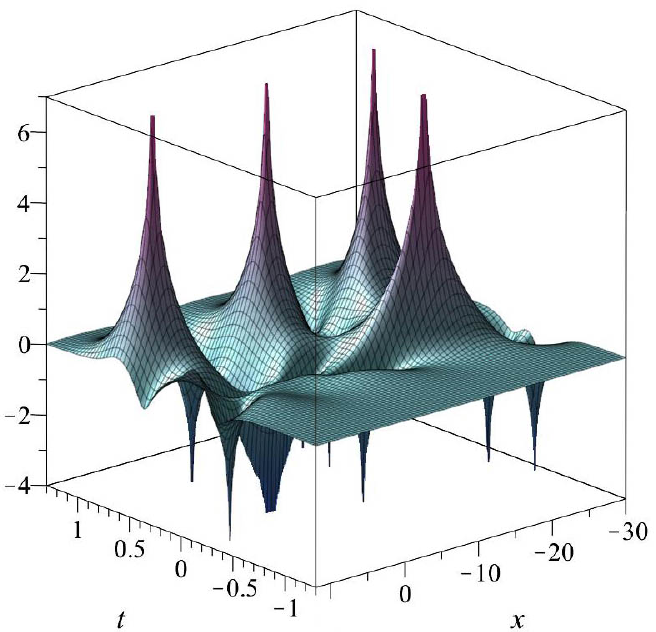}
\caption{$|\wt v|$ (left) and $\ln|\wt\rho |$ (right). \label{fig2}}
\end{figure}

{\it Case 2.}  When $p=1$, $C_1=\begin{bmatrix} c_{11}\\ c_{12} \end{bmatrix}$ and 
$C_2=\begin{bmatrix} 0 \\ c_{22} \end{bmatrix}$ $(c_{11},\, c_{12},\, c_{22}\in\BR)$,
our choice of non-diagonal $A$ leads to polynomials 
$$\g_1:=\I c_{12}x-2c_{11}-4c_{12}t/\bar a^2 \quad {\mathrm{and}} \quad
\g_2:=\I c_{12}x-2c_{11}-4c_{12}t/a^2$$
(in addition to the exponents)
in the formulas for $\wt v$ and $\wt \rho$. Namely, we have:
\begin{gather*}
\wt v=2a_1c_{22}\frac{
c_{22}^2(a_1\g_1-2c_{12})\E^{\I ax/2+2t/a}+c_{12}^2(a_1\g_2+2c_{12})\E^{-\I \bar ax/2-2t/\bar a}
}{
c_{22}^4\E^{\I a_1x+4a_1t/|a|^2}+c_{12}^4\E^{-\I a_1x-4a_1t/|a|^2}-2c_{12}^2c_{22}^2-a_1^2c_{22}^2\g_1\g_2
},
\\[1ex]
\wt\rho=\I +\frac {8\I a_1^2c_{22}^2}{|a|^4\big(c_{22}^4\E^{\I a_1x+4a_1t/|a|^2}+c_{12}^4\E^{-\I a_1x-4a_1t/|a|^2}-2c_{12}^2c_{22}^2-a_1^2c_{22}^2\g_1\g_2\big)^2}
\\ \qquad{} \times\Big(
 c_{22}^4(a_1\bar a\g_1+2\I a_2c_{12})(a_1a\g_2-2\I a_2c_{12})\E^{\I a_1x+4a_1t/|a|^2}\\ \qquad\quad{}
+c_{12}^4(a_1\bar a\g_1-2\I a_2c_{12})(a_1a\g_2+2\I a_2c_{12})\E^{-\I a_1x-4a_1t/|a|^2}\\ \qquad\quad{}
+2c_{12}^2c_{22}^2\big( -a_1^2(a_1^2-a_2^2)\g_1\g_2 + 32(\I c_{12}x-2c_{11})c_{12}a_2^2a_1^4|a|^{-4}t-4a_2^2c_{12}^2\big)
\Big).
\end{gather*}
The behaviour of  $\wt v$ and $\wt \rho$ is in this case more complicated, see Figure \ref{fig2},
where 
$$p=1, \quad a=\frac13+\frac15\I, \quad c_{11}=3, \quad c_{12}=1, \quad c_{22}=\frac12 .$$

\begin{figure}[h!]
\centering
\includegraphics[scale=1.0]{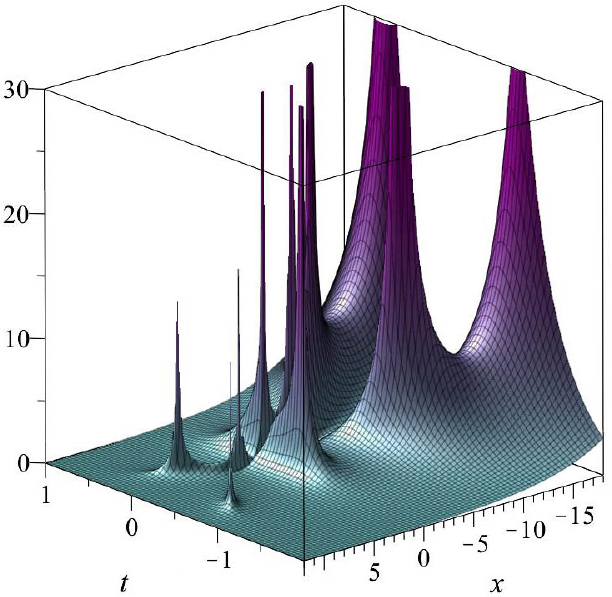}\qquad 
\includegraphics[scale=1.0]{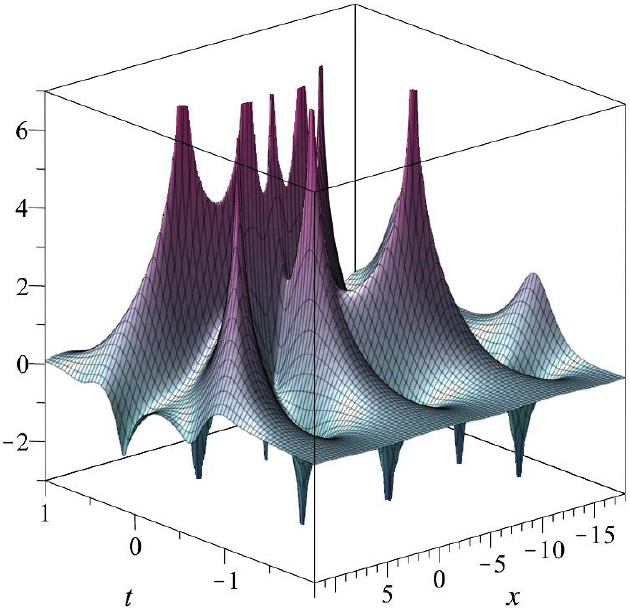}    
\caption{$|\wt v|$ (left) and $\ln|\wt\rho|$ (right). } \label{fig3}
\end{figure}

\begin{figure}[h!]
\centering
\includegraphics[scale=1.0]{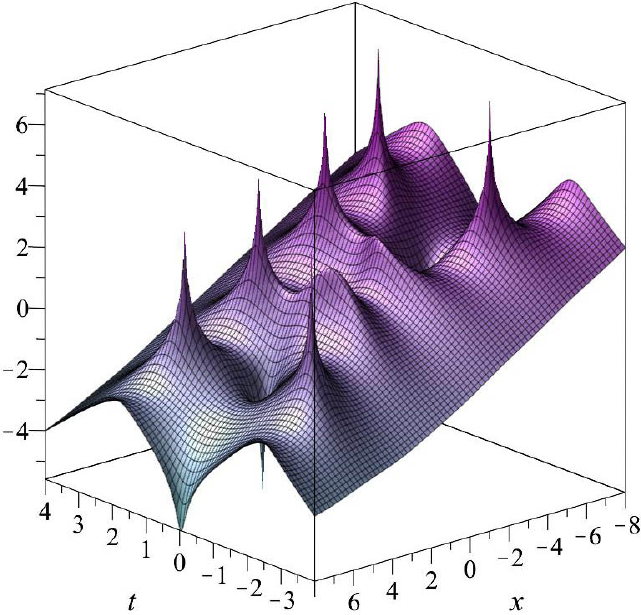}\quad\  
\includegraphics[scale=1.0]{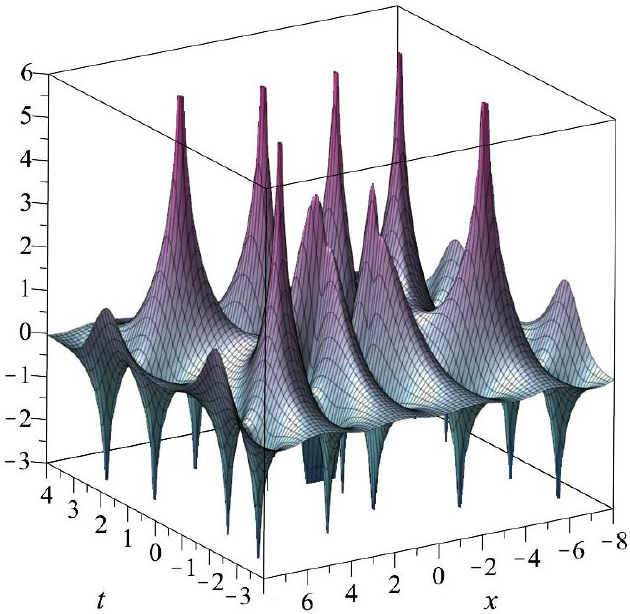}       
\caption{$\ln|\wt v|$ (left) and $\ln|\wt\rho|$ (right). }\label{fig4}
\end{figure}

\begin{figure}
\centering
\includegraphics[scale=1.0]{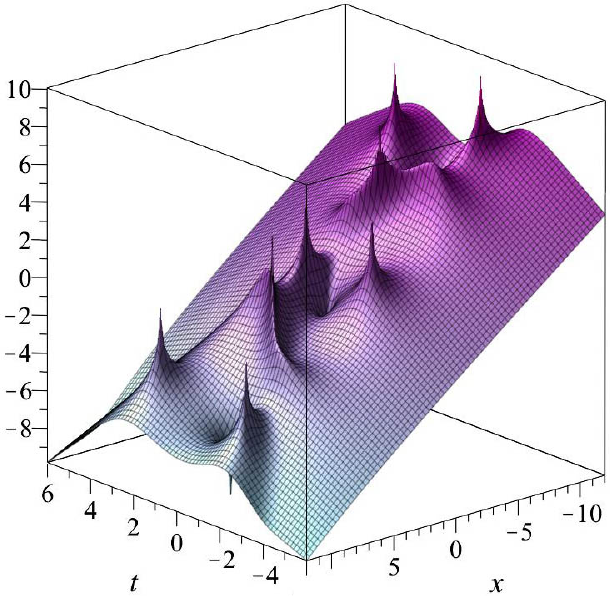}\qquad  
\includegraphics[scale=1.0]{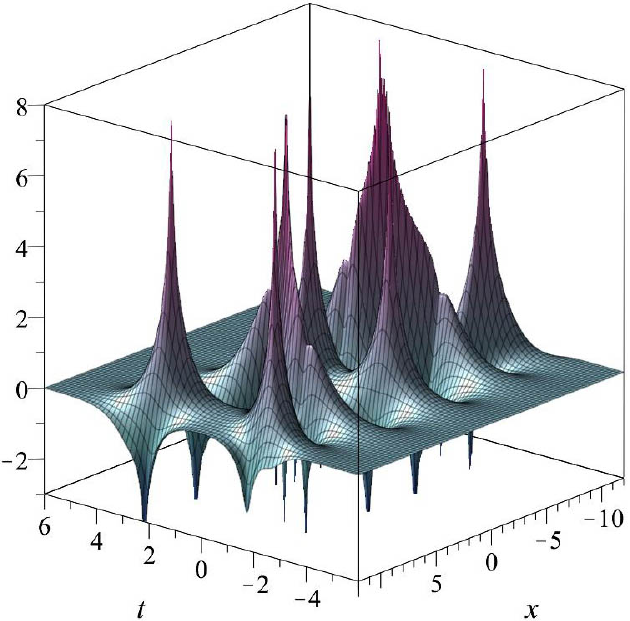}       
\caption{$\ln|\wt v|$ (left) and $\ln|\wt\rho|$ (right). }\label{fig5}
\end{figure}

In some other cases the formulas are more complicated and we restrict ourselves
to figures only. See Figure \ref{fig3}, where
$$p=1, \quad a=\frac12+\frac13\I , \quad C_1=\begin{bmatrix} 1 \\ 2\I \end{bmatrix}, \quad C_2=\begin{bmatrix} 3\I \\ 4 \end{bmatrix};$$
see Figure \ref{fig4}, where
$$p=0, \quad a=\frac32+\frac12\I , \quad C_1=\begin{bmatrix} 1+3\I \\3+ 2\I \end{bmatrix}, \quad C_2=\begin{bmatrix} 6+\I \\2-4 \I\end{bmatrix};$$
and see Figure \ref{fig5}, where
$$p=1, \quad a=1+\I , \quad C_1=\begin{bmatrix}\I \\ 1 \end{bmatrix}, \quad C_2=\begin{bmatrix} 1 \\  \I\end{bmatrix}.$$

\newpage

\bigskip

\noindent{\bf Acknowledgments.}  This research was supported by the
Austrian Science Fund (FWF) under Grant No. P29177. 

\newpage
 

\begin{flushright}
Faculty of Mathematics,
University
of
Vienna, \\
Oskar-Morgenstern-Platz 1, A-1090 Vienna,
Austria. \\
\vspace{1em}

R.O. Popovych, e-mail: {\tt roman.popovych@univie.ac.at}

\vspace{0.5em} 

A.L. Sakhnovich,
e-mail: {\tt oleksandr.sakhnovych@univie.ac.at}

\end{flushright}

\end{document}